\documentclass{amsart}
\usepackage{amssymb}
\usepackage{color}

\usepackage{a4}
\begin{document}
\def\nn{\nonumber}
\def\bea{\begin{array}}
\def\eea{\end{array}}
\def\beq{\begin{eqnarray}}
\def\eeq{\end{eqnarray}}
\newtheorem{theorem}{Theorem}[section]
\newtheorem{lemma}[theorem]{Lemma}
\newtheorem{remark}[theorem]{Remark}
\newtheorem{definition}[theorem]{Definition}
\newtheorem{corollary}[theorem]{Corollary}
\newtheorem{example}[theorem]{Example}
\makeatletter
\renewcommand{\theequation}{%
  \thesection.\alph{equation}}
\@addtoreset{equation}{section}
\makeatother
\title[Stability theorems for the bag]
{Stability theorems for chiral bag boundary conditions}
\author{P. Gilkey and K. Kirsten}
\begin{address}{PG: Math. Dept., University of Oregon, Eugene, Or 97403, USA}\end{address}
\begin{email}{gilkey@darkwing.uoregon.edu}\end{email}
\begin{address}{KK: Department of Mathematics, Baylor University \\
Waco, TX 76798, USA}\end{address}
\begin{email}{Klaus\_Kirsten@baylor.edu}\end{email}
\begin{abstract} We study asymptotic expansions of the smeared
$L^2$-traces $Fe^{-t P^2}$ and $FPe^{-tP^2}$,  where $P$ is an
operator of Dirac type and $F$ is an auxiliary smooth
endomorphism. We impose chiral bag boundary conditions depending
on an angle $\theta$. Studying the $\theta$-dependence of the
above trace invariants, $\theta$-independent pieces are
identified. The associated stability theorems allow one to show the
regularity of the eta function for the problem and to determine the most important heat
kernel coefficient on a four dimensional manifold.
\end{abstract}
\maketitle
\def\BB{{\mathcal{B}}}
\def\trl#1{\operatorname{Tr}_{L^2}\left\{#1\right\}}
\def\trv{\operatorname{Tr}_V}
\section{Introduction}\label{sect-1}
Let $M$ be a compact connected $m$-dimensional Riemannian manifold
with smooth boundary $\partial M$. Let $D$ be an operator of
Laplace type on a vector bundle $V$ over $M$. Let $D_\BB$ be the
realization of $D$ which is defined by a strongly elliptic
boundary condition $\BB$. Work of Greiner \cite{grei71-163} and
Seeley \cite{seel68-10-288,seel69-91-963} shows that the
fundamental solution of the heat equation $e^{-tD_\BB}$ is an
infinitely smoothing operator which is of trace class. Let $F\in
C^\infty(\operatorname{End}(V))$ be an auxiliary endomorphism of
$V$ which is used for localization. As $t\downarrow0$, there is a
complete asymptotic expansion
$$
\trl{Fe^{-tD_\BB}}\sim\sum_{n=0}^\infty a_n(F,D,\BB)t^{(m-n)/2}\,.
$$
The {\it heat asymptotics} $a_n$ are locally computable. There
exist suitable local endomorphisms $e_n(x,D)$ and
$e_{n,\nu}(y,D,\BB)$ of $V$ so that
$$
a_n(F,D,\BB)=\int_M\trv\left\{F(x)e_n(x,D)\right\}dx+
\sum_{\nu=0}^{n-1}\int_{\partial
M}\trv\left\{F^{(\nu)}e_{n,\nu}(y,D,\BB)\right\}dy\,.
$$
In this equation, $dx$ and $dy$ denote the Riemannian measures on
$M$ and on $\partial M$, respectively, and $F^{(\nu)}$ denotes the
$\nu^{\operatorname{th}}$ covariant derivative of $F$ with respect
to the inward unit normal using the canonical connection
determined by $D$.

These asymptotics can be computed quite explicitly; the principle
of ``not feeling the boundary'' shows that the interior invariants
$e_n(x,D)$ do not depend on the boundary condition chosen. The
invariants $e_n(x,D)$ vanish for $n$ odd and are known for
$n=0,2,4,6,8$, see for example \cite{Av90,Av91,AmBeOc89}. The
boundary invariants $e_{n,\nu}(y,D,\BB)$ are considerably more
subtle. The invariants for $n$ odd do not vanish. For mixed
boundary conditions they are known for $n=0,1,2,3,4,5$, and for
many other boundary conditions varying amounts of information are
known. We refer to \cite{gilk04,kirs01,vass03} for further details
concerning these formulas.

There are various stability theorems for these invariants that
play a crucial role in their evaluation. In Section \ref{sect-2},
we present results in the context of chiral bag boundary
conditions \cite{hras84-245-118,gold83}. To motivate these
stability results, we first give examples in the standard setting
which have proved to be crucial in past analysis.

The following stability formulae arise from conformal variations of the operator. Assertion (1)
follows from work of Branson and \O rsted \cite{BrOr86,BrOr88}; Assertion (2) was first observed
in \cite{BG90}.

\begin{theorem}\label{thm-1.1}
Let $\BB$ define an elliptic boundary condition for an operator
$D$ of Laplace type on a compact $m$-dimensional Riemannian
manifold $M$ with smooth boundary. Let $h\in C^\infty(M)$, let
$D_\varepsilon:=e^{-2\varepsilon h}D$, let $F\in
C^\infty(\operatorname{End}(V))$, and let
$F_\varepsilon:=e^{-2\varepsilon h}F$. Then
\begin{enumerate}
\item  $\partial_\varepsilon a_m(1,D_\varepsilon,\BB)=0$.
\item $\partial_\varepsilon a_{m-2}(F_\varepsilon,D_\varepsilon,\BB)=0$.
\end{enumerate}
\end{theorem}

\begin{proof} We proceed formally, the necessary analytic justification can be found in
\cite{gilk95b}:
\begin{eqnarray*}
&&\sum_{n=0}^\infty \partial_\varepsilon a_n(1,D_\varepsilon,\BB)t^{(n-m)/2}
\sim\partial_\varepsilon\trl{e^{-tD_{\varepsilon,\BB}}}\\
&=&-t\trl{(\partial_\varepsilon D_{\varepsilon ,\BB}
)e^{-tD_{\varepsilon,\BB}}}
=2t\trl{hD_{\varepsilon ,\BB} e^{-tD_{\varepsilon,\BB}}}\\
&=&-2t\partial_t\trl{he^{-tD_{\varepsilon,\BB}}}\sim
-2t\partial_t\sum_{n=0}^\infty  a_n(h,D_\varepsilon,\BB)t^{(n-m)/2}\\
&\sim&\sum_{n=0}^\infty (m-n)a_n(h,D_\varepsilon,\BB)t^{(n-m)/2}\,.
\end{eqnarray*}
Equating powers in the relevant asymptotic expansions yields
$$\partial_\varepsilon a_n(1,D_\varepsilon,\BB)=(m-n)a_n(h,D_\varepsilon,\BB)\,.$$
We set $n=m$ to complete the proof of Assertion (1).

Let $D_{\delta, \BB} :=D_\BB -\delta F$. We compute
\begin{eqnarray*}
&&\sum_{n=0}^\infty \partial_\delta a_n(1,D_\delta,\BB)t^{(n-m)/2}\sim
  \partial_\delta\trl{e^{-tD_{\delta,\BB}}}\\
&=&-t\trl{(\partial_\delta D_{\delta , \BB})e^{-tD_{\delta,\BB}}}=
  t\trl{Fe^{-tD_\delta,\BB}}\\
&\sim&t\sum_{n=0}^\infty a_n (F,D_\delta,\BB)t^{(n-m)/2}\,.
\end{eqnarray*}
Equating terms in the asymptotic expansions yields
\begin{equation}\label{eqn-1.a}
\partial_\delta a_n(1,D-\delta F,\BB)=a_{n-2}(F,D-\delta F,\BB)\,.
\end{equation}
We now consider the joint variation $D_{\delta,\varepsilon ,\BB
}:=e^{-2\varepsilon h}(D_\BB -\delta F)$. We use Assertion (1) and
Equation (\ref{eqn-1.a}) to see
\begin{eqnarray*}
&&\partial_\varepsilon a_m(1,e^{-2\varepsilon h}(D-\delta F),\BB)=0,\\
&&\partial_\delta a_m(1,e^{-2\varepsilon h}(D-\delta F),\BB)=a_{m-2}(e^{-2\varepsilon h}F,
e^{-2\varepsilon h}(D-\delta F),\BB)\,.
\end{eqnarray*}
Consequently:
\begin{eqnarray*}
0&=&\partial_\delta\partial_\varepsilon a_m(1,e^{-2\varepsilon h}(D-\delta F),\BB)
=\partial_\varepsilon\partial_\delta a_m(1,e^{-2\varepsilon h}(D-\delta F),\BB)\\
&=&\partial_\varepsilon a_{m-2}(e^{-2\varepsilon h}F,e^{-2\varepsilon h}(D-\delta F),\BB)\,.
\end{eqnarray*}
Assertion (2) follows by setting $\delta=0$. \end{proof}

Here is a brief guide to the remainder of this paper. In Section
\ref{sect-2}, we will use a similar strategy to establish
stability results for the zeta and eta invariants, when chiral bag
boundary conditions \cite{hras84-245-118,gold83}  are imposed. In
Section \ref{sect-3}, we will use these arguments to establish the
regularity at $s=0$ of the eta invariant. Although this result can
be derived from Theorem 2.3.5 \cite{GS83a}, it seemed worth giving
a straightforward and elementary argument adapted to the setting
at hand which is more conceptual in nature and which is of
interest in its own right; the discussion in \cite{GS83a} dealt
with a very general setting and was, perhaps, somewhat opaque. The
evaluation of invariants in arbitrary dimensions, even at the
level of $a_2$, has turned out to be extremely involved
\cite{EGK05}. As an application of the stability results, we will
see in Section \ref{sect-4} that it is relatively easy to find the
coefficient $a_4$ for the case $F=1$ in the most relevant case of
a four dimensional manifold. In Section \ref{sect-5}, we give a
perturbative result showing that the $a_4$ coefficient does
exhibit $\theta$ dependence in the general setting. In Section
\ref{sect-6}, we remove the assumption that $M$ is orientable and
in the final section, we make some concluding remarks.

\section{Chiral bag boundary conditions}\label{sect-2}

We start this section by describing the chiral bag boundary
conditions. Let $\gamma$ give the vector bundle $V$ a
$\operatorname{Clif}(TM)$ module structure. If $\{e_i\}$ is a
local orthonormal frame for the tangent bundle $TM$, then
$\gamma_i:=\gamma(e_i)$ forms a collection of skew-adjoint
matrices satisfying the Clifford commutation relations:
$$\gamma_i\gamma_j+\gamma_j\gamma_i=-2\delta_{ij}\operatorname{Id}_V\,.$$
Let $\nabla$ be a connection on $V$ and let $\psi$ be an auxiliary endomorphism of $V$. We form
the associated operator of Dirac type
$$P:= \gamma_i\nabla_{e_i}+\psi\,.$$
Assume that $m=2\bar m$ is even and that $M$ is oriented. Let
$$\tilde\gamma:=(\sqrt{-1})^{\bar
m}\gamma_1...\gamma_m$$ be the normalized orientation;
$\tilde\gamma$ is self-adjoint and
$\tilde\gamma^2=\operatorname{Id}_V$. It is the generalization of
$\gamma_5$ to arbitrary even dimension.

We assume as a compatibility condition that
$$\tilde\gamma P+P\tilde\gamma=0\,.$$
This means that if we decompose $V=V_+\oplus V_-$ into the $\pm1$ eigenbundles of $\tilde\gamma$,
then we may also decompose
$$
P=P_++P_-\quad\text{where}\quad P_\pm:C^\infty(V_\pm)\rightarrow V_\mp\,.
$$

Choosing suitable boundary conditions for $P$ is crucial. Whereas
spectral questions for many boundary conditions have been analyzed
in great detail \cite{gilk04,kirs01}, an understanding of
so-called chiral bag boundary conditions is still in its infancy;
see, however, \cite{EGK05,md03,wipf95-443-201}. These boundary
conditions are defined as follows. Set
$$\chi_\theta:= -\tilde\gamma e^{\theta\tilde\gamma}\gamma_m\quad\text{for}\quad
\theta\in\mathbb{R}\,.$$
Since $\gamma_m$ anti-commutes with $\tilde\gamma$, one has
$$\chi_\theta^2=\tilde\gamma
e^{\theta\tilde\gamma}\gamma_m\tilde\gamma
e^{\theta\tilde\gamma}\gamma_m
=\tilde\gamma\gamma_m\tilde\gamma\gamma_me^{-\theta\tilde\gamma}e^{\theta\tilde\gamma}
=\operatorname{Id}_V\,.
$$
We consider the projection operator on the $-1$ spectrum of
$\chi_\theta$, which is given by
$$\BB_\theta:=\textstyle\frac12(1-\chi_\theta)\,.$$
This defines suitable boundary conditions for $P$; let
$P_{\BB_\theta}$ be the realization of $P$. We set $D_{\BB_\theta}
=P_{\BB _\theta}^2$ and let the associated boundary operator be
$$
\BB_\theta^1\phi:=\BB_\theta\phi\oplus\BB_\theta P\phi\,.
$$
This is an elliptic boundary condition so the heat trace asymptotics are well defined
\cite{BGKS03}. To simplify the notation, we set
$$a_n(F,P,\BB_\theta):=a_n(F,D,\BB_\theta^1)\,.$$

The study of the coefficient $a_m$ is crucial to understanding the
chiral anomaly in the zeta function regularization; see, for
example, the discussion in \cite{md03}. The angle $\theta$
occurring in these boundary conditions is a substitute for
introducing small quark masses to drive the breaking of chiral
symmetry \cite{wipf95-443-201,duer97-255-333,duer99-273-1}. One of
the first papers where the chiral boundary conditions were
introduced is the work by Hrasko and Balog \cite{hras84-245-118},
with first applications to chiral bag models being presented in
\cite{gold83}.

The stability results we are going to prove naturally relate heat
equation invariants and eta invariants. The eta invariants measure
the spectral asymmetry of $P_{\BB_\theta}$; they are defined as
follows. If $F\in C^\infty(\operatorname{End}(V))$, then we may
expand
$$\trl{FP_{\BB_\theta}e^{-tP_{\BB_\theta}^2}}=\sum_{n=0}^\infty a_n^\eta(F,P,\BB_\theta)
t^{(n-m-1)/2}\,.$$ The eta invariants $a_n^\eta(F,P,\BB_\theta)$
are again locally determined. We refer to \cite{bene02-35-9343}
for a further discussion of the eta invariant in this setting.

We now come to one of the main results of this article.
\begin{theorem}\label{thm-2.1}
Let $P=\gamma_i\nabla_{e_i}+\psi$ be an operator of Dirac type on
a compact oriented smooth manifold of even dimension $m$ with
boundary. Let $F\in C^\infty(\operatorname{End}(V))$. Assume that
$P\tilde\gamma+\tilde\gamma P=0$ and that
$F\tilde\gamma+\tilde\gamma F=0$. Then
\begin{enumerate}
\item $\partial_\theta a_m(1,P,\BB_\theta)=0$.
\item $\partial_\theta a_m^\eta(1,P,\BB_\theta)=0$.
\item $\partial_\theta a_{m-1}^\eta(F,P,\BB_\theta)=0$.
\end{enumerate}
\end{theorem}

\begin{proof} The central technical point is to replace the given variation by one where the
boundary condition is held fixed, an idea that first occurred in
\cite{GS83}. We consider the gauge transformation defined by
$e^{\frac12\theta\tilde\gamma}$. Since $\tilde\gamma$
anti-commutes with $P$ and with $\gamma_m$, we may compute:
\begin{eqnarray*}
&&P_\theta:=e^{-\frac12\theta\tilde\gamma}Pe^{\frac12\theta\tilde\gamma}=e^{-\theta\tilde\gamma}
P,\\
&&e^{-\frac12\theta\tilde\gamma}\chi_\theta e^{\frac12\theta\tilde\gamma}
  = -e^{-\frac12\theta\tilde\gamma}\tilde\gamma e^{\theta\tilde\gamma}\gamma_m
   e^{\frac12\theta\tilde\gamma}=-e^{-\frac \theta 2\tilde\gamma}\tilde\gamma
e^{\frac \theta 2 \tilde\gamma}\gamma_m=\chi_0\,.
\end{eqnarray*}
Thus by gauge invariance,
$$a_n(1,P,\BB_\theta)=a_n(1,P_\theta,\BB_0)\,.$$

We have $\partial_\theta P_\theta=-\tilde\gamma P_\theta$; since we are multiplying on the left
by $\tilde\gamma$, the domain defined by the boundary condition $\BB_0$ is not perturbed. We can
now argue exactly as in the proof of Theorem \ref{thm-1.1} to see:
\begin{eqnarray*}
&&\partial_\theta\trl{e^{-tP^2_{\theta,\BB_0}}}
=-2t\trl{\partial_\theta(P_{\theta,\BB_0})P_{\theta,\BB_0}e^{-tP_{\theta,\BB_0}^2}}\\
&=&2t\trl{\tilde\gamma P^2_{\theta,\BB_0}e^{-tP_{\theta,\BB_0}^2}}
=-2t\partial_t\trl{\tilde\gamma e^{-tP_{\theta,\BB_0}^2}}\,.
\end{eqnarray*}
Expanding in an asymptotic series then yields:
\begin{eqnarray*}
&&\sum_{n=0}^\infty \partial_\theta a_n(1,P_\theta,\BB_0)
t^{(n-m)/2}\sim
\sum_{n=0}^\infty-2t\partial_t\left\{a_n(\tilde\gamma,P_\theta,B_0)t^{(n-m)/2}\right\}\\
&\sim&
\sum_{n=0}^\infty (m-n)a_n(\tilde\gamma,P_\theta,\BB_0)t^{(n-m)/2}\,.
\end{eqnarray*}
Equatating terms in the relevant asymptotic expansions yields:
\begin{equation}\label{eqn-2.a}
\partial_\theta a_n(1,P,\BB_\theta)=\partial_\theta a_n(1,P_\theta,\BB_0)
=(m-n)a_n(\tilde\gamma,P_\theta,\BB_0)\,.
\end{equation}
Assertion (1) then follows by setting $n=m$ in Equation (\ref{eqn-2.a}).

The argument is similar to prove Assertion (2). We compute:
\begin{eqnarray*}
&&\sum_{n=0}^\infty \partial_\theta a_n^\eta(1,P,\BB_\theta)=
  \sum_{n=0}^\infty \partial_\theta a_n^\eta(1,P_\theta,\BB_0)\\
&\sim&\partial_\theta\trl{P_{\theta,\BB_0}e^{-tP_{\theta,\BB_0}^2}}
=\trl{\partial_\theta(P_{\theta,\BB_0})(1-2tP_{\theta,\BB_0}^2)e^{-tP_{\theta,\BB_0}^2}}\\
&=&(1+2t\partial_t)\trl{-\tilde\gamma P_{\theta,\BB_0}e^{-tP_{\theta,\BB_0}^2}}\\
&\sim&-(1+2t\partial_t)\sum_{n=0}^\infty a_n^\eta(\tilde\gamma,P_\theta,\BB_0)^{(n-m-1)/2}\\
&\sim&\sum_{n=0}^\infty(m-n)a_n^\eta(\tilde\gamma,P_\theta,\BB_0)t^{(n-m-1)/2}
.
\end{eqnarray*}
This yields the relation
$$
\partial_\theta a_n^\eta(1,P,\BB_\theta)=(m-n)a_n^\eta(\tilde\gamma,P_\theta,\BB_0)\,.
$$
Assertion (2) follows by setting $n=m$.

To prove Assertion (3), consider a variation of the form $P_\delta:=P-\delta F$. We argue as in
the proof of Theorem
\ref{thm-1.1} (2) to compute:
\begin{eqnarray*}
&&\sum_{n=0}^\infty \partial_\delta
a_n(1,P_\delta,\BB_\theta)t^{(n-m)/2}\sim\partial_\delta\trl{ e^{-tP_{\delta,\BB_\theta}^2}}\\
&=&-2t\trl{(\partial_\delta
P_{\delta,\BB_\theta})P_{\delta,\BB_\theta}e^{-tP_{\delta,\BB_\theta}^2}}
=2t\trl{FP_{\delta,\BB_\theta}e^{-tP_{\delta,\BB_\theta}^2}}\\
&\sim&2\sum_{k=0}^\infty a_k^\eta(F,P,\BB_\theta)t^{(k-m+1)/2}\,.
\end{eqnarray*}
Equating terms in the asymptotic expansion then yields
\begin{equation}\label{eqn-2.b}
2a_k^\eta(F,P,\BB_\theta)=\partial_\delta a_{k+1}(1,P_\delta,\BB_\theta)\,.
\end{equation}
Since $P_\delta\tilde\gamma+\tilde\gamma P_\delta=0$, we can apply
Theorem \ref{thm-2.1} (1). Differentiating Equation
(\ref{eqn-2.b}) with respect to $\theta$, setting $k+1=m$ and
$\delta=0$, then yields the desired result.\end{proof}

\begin{remark}\label{rmk-2.2}\rm It is natural to try a similar approach to study $\partial_\theta
a_{m-1}(F,P,\BB_\theta)$. In fact, however, this fails. One would compute:
\begin{eqnarray*}
&&\sum_{n=0}^\infty\partial_\delta a_n^\eta(1,P_\delta,\BB_\theta)t^{(n-m-1)/2}\sim
\partial_\delta\trl{P_{\delta,\BB_\theta}e^{-tP_{\delta,\BB_\theta}^2}}\\
&=&\trl{(\partial_\delta P_{\delta,\BB_\theta})(1-2tP_{\delta,\BB_\theta}^2)
   e^{-tP_{\delta,\BB_\theta}^2}}\\
&=&-(1+2t\partial_t)\trl{F
   e^{-tP_{\delta,\BB_\theta}^2}}\\
&\sim&-(1+2t\partial_t)\sum_{k=0}^\infty
    a_k(F,P_\delta,\BB_\theta)t^{(k-m)/2}\\
&\sim&\sum_{k=0}^\infty
(m-k-1)a_k(F,P_\delta,\BB_\theta)t^{(k-m)/2}\,.
\end{eqnarray*}
Equating terms in the asymptotic expansions then yields
\begin{equation}\label{eqn-2.c}
\partial_\delta a_{k+1}^\eta(1,P_\delta,\BB_\theta)=(m-k-1)a_k(F,P_\delta,\BB_\theta)\,.
\end{equation}
Setting $k+1=m$ yields no information about $a_{m-1}(F,P_\delta,\BB_\theta)$.
\end{remark}

\begin{remark} \rm Assertion (1) of Theorem \ref{thm-2.1} explains
certain observations made in the literature. For example, the
calculation of the heat kernel coefficients on the ball with $F=1$
showed, in general dimensions, a strong dependence on $\theta$. It
was noticed that in $m=2$, respectively $m=4$, this dependence
disappeared when $a_2(1,P,\BB_\theta)$, respectively $a_4
(1,P,\BB_\theta )$, was considered \cite{giam02}. The above result
shows why this must be the case.

The same pattern could be observed in the coefficient $a_2
(1,P,\BB_\theta )$ found in \cite{EGK05}. Whereas in general the
$\theta$-dependence enters in terms of hypergeometric functions,
these terms are multiplied by $(m-2)$ and they thus vanish at
$m=2$; see Equation (2d) in \cite{EGK05}.
\end{remark}

\setcounter{section}{2}\section{Regularity of the eta invariant at
s=0}\label{sect-3} In this section we apply Theorem \ref{thm-2.1}
(2) to discuss the regularity of the eta invariant at $s=0$.

The eta invariant was introduced by Atiyah et. al. \cite{ABS73}
for closed manifolds. We extend those definitions to the setting
at hand as follows. Let
$$\eta(s,P,\BB_\theta):=\trl{P_{\BB_\theta} (P_{\BB_\theta}^2)^{-(s+1)/2}}\,.$$
(There is a bit of fuss with the $0$-spectrum that can be dealt
with appropriately). The standard pseudo-differential calculus,
see, for example, the discussion in \cite{GS83}, can be used to
show $\eta(s,P,\BB_\theta)$ has a meromorphic extension to
$\mathbb{C}$ with isolated simple poles with locally computable
residues. If $P_{\BB_\theta}$ is self-adjoint with respect to some
Hermitian metric on $V$ and if $\{\lambda_\nu\}$ is the spectrum
of $P_{\BB_\theta}$, each eigenvalue being repeated according to
multiplicity, then
$$
\eta(s,P,\BB_\theta):=\sum_{\lambda_\nu\ne0}\operatorname{sign}(\lambda_\nu)|\lambda_\nu|^{-s}
\,.$$ A-priori, $s=0$ need not be a regular value of $\eta$.
However, the pole of the eta function at $s=0$ can be related to
the trace invariant $a_m^\eta(1,P,\BB_\theta)$:
$$\operatorname{Res}_{s=0}\eta(s,P , \BB_\theta)=
\textstyle2\Gamma(\frac12)^{-1}a_m^\eta(1,P,\BB_\theta)\,.$$ We
will show that $a_m^\eta(1,P,\BB_\theta)=0$ and thus $\eta$ is
regular at $s=0$. One then sets
$$
  \eta(P, \BB_\theta ):=\textstyle\frac12\big\{\eta(s,P,\BB_\theta)
+\dim\ker\{P_{\BB_\theta}\}\big\}\big|_{s=0} \in
\mathbb{C}/\mathbb{Z}
$$
as a global measure of the spectral asymmetry of $P$. For closed manifolds, this
invariant plays a central role in the index theorem for manifolds
with boundary \cite{ABS73}. It can also be used to study
$K$-theory groups and equivariant bordism groups \cite{G99} and to
study metrics of positive scalar curvature \cite{BGS97}. Thus it
is important to understand this invariant in the context of bag
boundary conditions.

\begin{theorem}\label{thm-3.1}
Let $P$ be an operator of Dirac type on a compact oriented smooth
manifold of even dimension $m$ with boundary. Then
$a_m^\eta(1,P,\BB_\theta)=0$.
\end{theorem}

\begin{proof}  Let $Q$ be an auxilary first order
partial differential operator on $V$. There is a complete asymptotic
expansion, see for example Lemma 2.6 \cite{GS83},
$$\trl{Qe^{-tP_{\BB_\theta}^2}}=\sum_{n=0}^\infty a_n^\eta(Q,P,\BB_\theta)
t^{(n-m-1)/2}\,.$$ Let $P_\varepsilon$ be a smooth $1$-parameter
family of operators of Dirac type. We assume the leading symbol of
$P_\varepsilon$ agrees with the leading symbol of $P$ on $\partial
M$. Thus $\BB_\theta$ defines an elliptic boundary condition for
this family as well. Using the notation $\dot P_\varepsilon =
\partial_\varepsilon (P_{\varepsilon, \BB_\theta})$, we compute:
\begin{eqnarray*}
&&\sum_{n=0}^\infty\partial_\varepsilon
a_n^\eta(1,P_\varepsilon,\BB_\theta)t^{(n-m-1)/2}
\sim\partial_\varepsilon\trl{P_{\varepsilon,\BB_\theta}e^{-tP_{\varepsilon,\BB_\theta}^2}}\\
&=&\trl{ \dot P _\varepsilon (1-2tP_{\varepsilon,\BB_\theta}^2)
e^{-tP_{\varepsilon,\BB_\theta}^2}}\\
&=&(1+2t\partial_t)\trl{\dot P_{\varepsilon}e^{-tP_{\varepsilon,\BB_\theta}^2}}\\
&\sim&\sum_{n=0}^\infty(1+2t\partial_t)a_n^\eta(\dot
P_\varepsilon,P_\varepsilon,\BB_\theta)
  t^{(n-m-1)/2}\\
&\sim&\sum_{n=0}^\infty(n-m)a_n^\eta(\dot
P_\varepsilon,P_\varepsilon,\BB_\theta)
  t^{(n-m-1)/2}\,.
\end{eqnarray*}
Equating coefficients in the asymptotic expansion then yields
$$\partial_\varepsilon a_n^\eta(1,P_\varepsilon,\BB_\theta)=(n-m)a_n^\eta(\dot P_\varepsilon,
P_\varepsilon,\BB_\theta)\,.$$ Setting $n=m$ then yields the well
known variational principle:
$$
\partial_\varepsilon a_m^\eta(1,P_\varepsilon,\BB_\theta)=0\,.
$$
Equation (\ref{eqn-2.c}) is a special case of this construction.

Using a partition of unity, we can construct a smooth
$1$-parameter family of operators $P_\varepsilon$ which are of
Dirac type with respect to Riemannian metrics $g_\varepsilon$ so
that $P_0=P$, so that $P_\varepsilon=P$ on $\partial M$, and so
that the metric defined by $P_1$ is product near the boundary.
Since $a_m^\eta$ is unchanged by this variation, we may assume
without loss of generality that the Riemannian metric in question
is product near the boundary in proving Theorem \ref{thm-3.1}.

We say that a connection $\nabla$ on $V$ is compatible if $\nabla\gamma=0$; such
connections always exist \cite{BG92}.  Choose such a connection on $V$ and
expand
$P=\gamma_i\nabla_{e_i}+\psi$. Let
$$P_\varepsilon:=\gamma_i \nabla_{e_i}+\varepsilon\psi\,.$$
Since $a_m^\eta(1,P_\varepsilon,\BB_\theta)$ is independent of
$\varepsilon$, setting $\varepsilon=0$ shows we may assume
$\psi=0$ in proving Theorem \ref{thm-3.1}. Since $\psi=0$,
$\tilde\gamma$ anti-commutes with $P$. Thus Theorem 2.1 (2)
permits to restrict to the case $\theta=0$.

By results in \cite{GS83}, there are local invariants
$a_n^\eta(x,P)$ and $a_n^\eta(y,P,\BB_0)$ so
$$a_n^\eta(1,P,\BB_\theta)=\int_Ma_n^\eta(x,P)dx+\int_{\partial M}a_n^\eta(y,P,\BB_0)dy\,.$$
The interior invariant $a_n^\eta (x,P)$ is homogeneous of order
$n$ in the jets of the symbol of $P$  and the boundary invariant
$a_n^\eta(y,P,\BB_0)$ is homogeneous of order $n-1$ in the jets of
the symbol of $P$. There is a parity constraint, by which the
interior invariants $a_n^\eta (x,P)$ vanish if $n$ is {\it even}.
In particular, this invariant plays no role in the study of
$a_m^\eta(1,P,\BB_0)$.

Let $U$ be a small contractible open neighborhood of $y\in\partial
M$. Since $\nabla$ is a compatible connection, Theorem 1.4
\cite{BG92} shows one may decompose
\begin{equation}\label{eqn-3.b}
V|_U=\Delta\otimes V_1\quad\text{and}\quad
\nabla|_U=\nabla^s\otimes\operatorname{Id}+\operatorname{Id}\otimes\nabla^1
\end{equation}
where $\nabla^s$ is the spin connection on the spin bundle
$\Delta$ and where $\nabla^1$ is an auxiliary connection on an
auxiliary vector bundle $V_1$. Since the structures are product
near the boundary, the boundary invariant $a_m^\eta(y,P,\BB_0)$ is
a universal polynomial in the covariant derivatives of the
curvature tensors of $(\nabla^s,\nabla^1)$; the invariants in
question are local so the possible failure of the decomposition
given in Equation (\ref{eqn-3.b}) to exist globally plays no role.
Since $\tilde\gamma$ depends on the orientation, the structure group is $SO(m-1)$. The
normal vector plays no role since the structures are product near the boundary; these
invariants are defined by the structures on $\partial M$.

The invariant $a_m^\eta(y,P,B_0)$ is homogeneous of odd degree $m-1$ in
the derivatives of the metric on the boundary and in the derivatives of
the connection $1$-form of $\nabla^1$. We can decompose
$a_m^\eta=a_{m,+}^\eta+a_{m,-}^\eta$ where $a_{m,+}^\eta$ is an $O(m-1)$
invariant and where $a_{m,-}^\eta$ changes sign if the orientation is reversed.
By H. Weyl's theorem, all local scalar $O(m-1)$ invariants are constructed by
taking traces and contracting tangential indices in pairs. Since
$a_{m,+}^\eta(y,P,\BB_0)$ is homogeneous of order $m-1$ and since $m-1$ is odd,
this local invariant must vanish; there are no odd order invariants in the jets
of the curvature of the metric $g$ and of the auxiliary connection $\nabla^1$ since the
structures are product near the boundary.

The invariant $a_{m,-}^\eta$ changes sign if the orientation is reversed. The
invariance theory used in the heat equation proof of the local index theorem
shows that this invariant must vanish as well, see for example the discussion
in \cite{gilk95b,G75}.\end{proof}

\begin{remark}\label{rmk-3.2}
\rm It is crucial in the analysis performed above that invariants
involving the endomorphism $\psi$ or the second fundamental form
$L$ are eliminated by performing a relevant homotopy. For example,
if $m=4$, the invariants $L_{aa}R_{ijij}$ and
$\operatorname{Tr}(\psi_{;m}\psi)$ are $O(3)$ invariants on the
boundary which could a-priori enter. Let $\varepsilon$ be the
totally anti-symmetric tensor. The invariant
$\{\varepsilon_{abc}L_{ad}R_{bcdm}\}$ is an $SO(3)$ invariant
which changes sign if the orientation is reversed; this invariant
plays a crucial role studying the axial anomaly for a Euclidean
Taub-NUT metric, see, for example, the discussion in \cite{EGH78}.
In our setting, this invariant is eliminated from consideration
since we perturb the problem so that the second fundamental form
vanishes.
\end{remark}

\begin{remark}\label{rmk-3.3} \rm We have shown that
$\eta(s,P,\BB_\theta)$ is regular at $s=0$. We defined
$$\eta(P,\BB_\theta):=\textstyle\frac12\big\{\eta(s,P,\BB_\theta)
+\dim\ker\{P_{\BB_\theta}\}\big\}\big|_{s=0}
\in\mathbb{C}/\mathbb{Z}
$$
as a measure of the global spectral asymmetry of $P$. Although
this invariant is not locally computable, the derivative is
locally computable. If $P_\varepsilon:=P+\varepsilon\Xi$ where
$\Xi$ is a $0^{\operatorname{th}}$ order operator, then one has,
see for example \cite{gilk95b}, that
$$\partial_\varepsilon \eta(P_\varepsilon,\BB_\theta) |_{\epsilon =0} =
2\Gamma({\textstyle\frac12})^{-1}a_{m-1}(\Xi,P,\BB_\theta)\,.$$ If
$\partial M=\emptyset$, then $a_{m-1}(\Xi,P,\BB_\theta)=0$ and
$\eta(P_\varepsilon,\BB_\theta)$ is independent of $\varepsilon$.
This is not, however, the case if $\partial M$ is non-empty. Take
$\Xi:=-f_o\tilde\gamma\gamma_m$ where $f_o$ vanishes to second
order on the boundary and where $f_o$ is supported near the
boundary. Set $\theta=0$, and set $m=4$. By Theorem IV.1
\cite{BGV98},
$$
a_3(\Xi,P,\BB_0)=(4\pi)^{-3/2}(384)^{-1}\int_{\partial M}30\dim(V)f_{o;mmm}dy\,.
$$
Thus the eta invariant is {\bf not} a homotopy invariant on an even dimensional manifold
 with these boundary conditions
if the boundary is non-empty.
\end{remark}

\section{Evaluating the coefficient $a_4 (1,P,\BB_\theta)$}\label{sect-4}

As a further immediate application of Theorem \ref{thm-2.1} we are
able to evaluate $a_4 (1,P,\BB_\theta )$ for $P=\gamma_i
\nabla_{e_i}+\psi$ when $m=4$ and when $P\tilde\gamma+\tilde\gamma
P=0$; this is the most relevant dimension. Given that
$\partial_\theta a_4 (1,P,\BB_\theta )=0$ for $m=4$, we have $a_4
(1,P,\BB_\theta ) = a_4 (1,P,\BB_0)$. However, for $\theta =0$ the
boundary conditions reduce to standard boundary conditions of
mixed type. To use the known results for mixed boundary conditions
\cite{gilk04,kirs01} we introduce some notation. We adopt the
Bochner formalism. If $D$ is an operator of Laplace type, there is
a unique connection $\nabla^D$ and a unique endomorphism $E$ so
that
$$D=-g^{ij} \nabla_i^D \nabla_j^D -E\,.$$
Let $\chi$ be an auxiliary Hermitian endomorphism used to define
the splitting of the bundle $V$. With the projectors
$$\Pi_\pm =\textstyle \frac 1 2 \left( 1\pm \chi \right) ,$$ the mixed
boundary operator is defined as \beq \BB_0^1 \phi = \Pi _- \phi
\left| _{\partial M} \oplus ( \nabla_m^D + S) \Pi _+ \phi \right|
_{\partial M} =0.\nn \eeq Let $\Omega_{ij}^D$ be the curvature of
the connection $\nabla^D$. Furthermore, let $\tau := R_{ijji}$ be
the scalar curvature, with the convention that the components
$R_{ijkl}$ of the curvature of the Levi-Civita connection are such
that $R_{1212}=-1$ for the standard metric on $S^2$. Finally, let
$\rho_{ij}=R_{ikkj}$ be the Ricci tensor and $L_{ab}$ the second
fundamental form and let `$;$' (resp. `$:$') denote the covariant
derivative (resp. tangential covariant derivative) with respect to
the connection $\nabla^D$ and the Levi-Civita connection of $M$
(resp. of $\partial M$). Using this notation, furthermore $\rho ^2
= \rho_{ij} \rho_{ij}$ and $R^2 = R_{ijkl} R_{ijkl}$, the
following result has been shown in \cite{BG90,dim95}.
\begin{theorem}\label{thm-4.1}
\begin{eqnarray*} &&a_4(1,D,\BB_0^1 ) \\&=&\frac1{ 360(4\pi)^{m/2}}
\int_M\operatorname{Tr}_V[60E_{;kk}+60\tau E+180E^2+
30\Omega_{ij}^D\Omega_{ij}^D+12\tau_{;kk}\\&&+5\tau^2-2\rho^2+2R^2]dx\\
&&+\int\limits_{\partial M} \operatorname{Tr}_V   \left[ \left( 240
\Pi _+ - 120 \Pi _- \right) E_{;m} + \left( 42 \Pi _+ - 18 \Pi _-
\right)
\tau_{;m} \right.\nn\\
& &+120 E L_{aa} + 20 \tau L_{aa} - 4 \rho_{mm} L_{aa} - 12
R_{ambm} L_{ab} + 4 R_{abcb}L_{ac} \nn\\
& &+\frac 1 {21} \left( 280 \Pi _+ + 40 \Pi _-\right) L_{aa}
L_{bb} L_{cc} + \frac 1 {21} \left( 168 \Pi _+ - 264 \Pi _-\right)
L_{ab} L_{ab} L_{cc} \nn\\
& &+\frac 1 {21} \left( 224 \Pi _+ + 320 \Pi _- \right) L_{ab}
L_{cb} L_{ac} + 720 SE + 120 S\tau \nn\\
& &+144 S L_{aa} L_{bb} + 48 S L_{ab} L_{ab} + 480 S^2 L_{aa} +
480
S^3\nn\\
& &\left. +60 \chi \chi _{:a} \Omega_{am}^D - 12 \chi_{:a}
\chi_{:a} L_{bb} - 24 \chi_{:a} \chi_{:b} L_{ab} -120 \chi_{:a}
\chi_{:a} S \right]dy\,\, .\nn\end{eqnarray*}
\end{theorem}
For the present case under investigation, we have by
\cite{BG92,BGG92} that writing $$P^2 = (\gamma_i \nabla_i + \psi
)^2 := -g_{ij} \nabla_i ^D \nabla_j^D - E , $$ the induced
connection $\nabla_i^D$ satisfies
$$\nabla_i^D = \nabla_i +
\omega_i\quad \text{where}\quad\omega_i = - \textstyle\frac 1 2 (\psi
\gamma_i + \gamma_i \psi ) .$$ Furthermore, in terms of the
induced connection $\nabla_i^D$, $P$ can be written as
$$P=\gamma_i \nabla_i^D + \phi \quad\text{where}\quad\phi = \psi -
\gamma_i \omega_i.$$ To state all ingredients needed for the
actual evaluation of Theorem 4.1 for the chiral bag boundary
condition with $\theta =0$, let $\Omega_{ij}$ be the curvature of
the compatible connection. Then, for any concrete example, the
formula in Theorem \ref{thm-4.1} can be evaluated from
\cite{BG92,BGG92}
\begin{eqnarray*}
\Omega_{ij}^D &=& \Omega_{ij} + \nabla_i \omega_j - \nabla_j
\omega_i + \omega_i
\omega_j - \omega_j \omega_i, \nn\\
E&=& -\textstyle\frac 1 2 \gamma_i \gamma_j \Omega_{ij} ^D - \gamma_i
\nabla_i \phi - \phi ^2 , \nn\\
S&=& \textstyle\frac 1 2 \Pi _+ (-\gamma_m \psi + \psi \gamma_m - L_{aa})
\Pi_+ , \nn\\
\chi &=& -\tilde \gamma \gamma_m , \quad \quad \chi_{:a} = \tilde
\gamma \gamma_b L_{ab} + \omega_a \chi - \chi \omega_a\,.
\end{eqnarray*}
\section{A perturbative result -- $\partial_\theta a_m(1,P,\BB_\theta)|_{\theta=0}$}\label{sect-5}
In this section, we drop the assumption
$P\tilde\gamma+\tilde\gamma P=0$ and examine, at least in part,
the $\theta$ dependence.

\begin{theorem}\label{thm-5.1}Expand $P=\gamma_i\nabla_{e_i}+\psi$ where $\nabla$ is compatible.
Let $\psi=\psi_o+\psi_e$
where
$\psi_o$ anti-commutes with
$\tilde\gamma$ and $\psi_e$ commutes with $\tilde\gamma$. Then:
\begin{enumerate}
\item $\partial_\theta\{a_m(1,P,\BB_\theta)\}|_{\theta=0}=
-2a_{m-1}^\eta(\tilde\gamma\psi_e,P,\BB_0)$.
\item Let $m=4$, let $d=\dim(V)$, let $\psi_o=f_0\tilde\gamma\gamma_m$,
and let $\psi_e=f_e\tilde\gamma$ where $f_o$ and $f_e$ are smooth
functions on $M$ and where $f_o$ is supported near $\partial M$.
Then:
$$\partial_\theta\{a_m(1,P,\BB_\theta)\}|_{\theta=0}=\frac d{16\pi^2}
\int_{\partial
M}\{2f_ef_oL_{aa}+2f_e^2f_o-2f_{e;m}f_o-f_ef_{o;m}\}dy\,.$$
\end{enumerate}\end{theorem}

\begin{remark}\rm
Assertion (2) shows that $a_m(1,P,\BB_\theta)$ exhibits non-trivial $\theta$ dependence in general.
\end{remark}
\begin{proof}
 Consider the gauge transformation defined by
$e^{\frac12\theta\tilde\gamma}$:
\begin{eqnarray*}
&&P_\theta:=e^{-\frac12\theta\tilde\gamma}(P-\psi_e+\psi_e)e^{\frac12\theta\tilde\gamma}=e^{-\theta\tilde\gamma}
(P-\psi_e)+\psi_e\\
&&\qquad=e^{-\theta\tilde\gamma}P+(1-e^{-\theta\tilde\gamma})\psi_e,\\
&&e^{-\frac12\theta\tilde\gamma}\chi_\theta e^{\frac12\theta\tilde\gamma}
  = -e^{-\frac12\theta\tilde\gamma}\tilde\gamma e^{\theta\tilde\gamma}\gamma_m
   e^{\frac12\theta\tilde\gamma}= -e^{-\frac \theta 2\tilde\gamma}\tilde\gamma
e^{\frac \theta 2 \tilde\gamma}\gamma_m=\chi_0\,.
\end{eqnarray*}
By gauge invariance and a direct calculation, we have
$$a_n(1,P,\BB_\theta)=a_n(1,P_\theta,\BB_0)\quad\text{and}\quad
\partial_\theta P_\theta|_{\theta=0}=-\tilde\gamma P+\tilde\gamma\psi_e\,.$$
Since we are multiplying
on the left by $\tilde\gamma$, the domain of
$P_\theta$ defined by the boundary condition
$\BB_0$ is not perturbed. We compute:
\begin{eqnarray*}
&&\partial_\theta\trl{e^{-tP^2_{\theta,\BB_0}}}\bigg|_{\theta=0}
=-2t\trl{\partial_\theta(P_{\theta,\BB_0})P_{\BB_0}e^{-tP_{\BB_0}^2}}\bigg|_{\theta=0}\\
&=&-2t\trl{(-\tilde\gamma P^2_{\BB_0}+\tilde\gamma
\psi_eP_{\BB_0})e^{-tP_{\BB_0}^2}} \\
&=&-2t\partial_t\trl{\tilde\gamma
e^{-tP_{\BB_0}^2}}
  -2t\trl{\tilde\gamma\psi_e P_{\BB_0}e^{-tP_{\BB_0}^2}}
\,.
\end{eqnarray*}

Expanding in an asymptotic series then yields:
\begin{eqnarray*}
&&\sum_{n=0}^\infty\partial_\theta a_n(1,P_\theta,\BB_0) t^{(n-m)/2}\bigg|_{\theta=0}\\
&\sim&
-2t\partial_t\sum_{j=0}^\infty
a_j(\tilde\gamma,P,\BB_0)t^{(j-m)/2}-2t\sum_{k=0}^\infty
a_k^\eta(\tilde\gamma\psi_e,P,\BB_0)t^{(k-m-1)/2}\,.
\end{eqnarray*}
Equating terms in the relevant asymptotic expansions yields:
\begin{eqnarray*}
&&\partial_\theta a_n(1,P,\BB_\theta)\big|_{\theta=0}=\partial_\theta
a_n(1,P_\theta,\BB_0)\big|_{\theta=0}\\
&=&(m-n)a_n(\tilde\gamma,P,\BB_0)-2a_{n-1}^\eta(\tilde\gamma\psi_e,P,\BB_0)\,.
\end{eqnarray*}
Assertion (1) then follows by setting $n=m$.

We specialize this to the case $m=4$. We apply Theorem 12 \cite{BGG92}. Set
$$\begin{array}{lll}
\psi:=f_o\tilde\gamma\gamma_m+f_e\tilde\gamma,&\psi_e:=f_e\tilde\gamma,&\psi_o:=f_o\tilde\gamma\gamma_m ,\\
\chi:=-\tilde\gamma\gamma_m,&
\Psi:=\gamma_i\psi\gamma_i,&
\Psi_T:=\gamma_a\psi\gamma_a\,.
\end{array}$$
After adjusting for the different sign convention employed,
setting $m=4$, and using $f=\tilde \gamma \psi_e =f_e$ scalar, we
find \begin{eqnarray*} &&
a_3^\eta(f,P,\BB_0)=-\frac1{192\pi^2}\bigg\{\int_Mf_e\trv
\big\{\{6\psi_{;i}+3\Psi\gamma_i\psi\}_{;i}\\
&&\qquad+\{-\tau  \psi-6 \gamma  _{ i} \gamma  _{ j}W _{ ij} \psi +6 \psi  \psi  _{ ;i}
\gamma  _{ i}
  -3 \psi  \psi  \Psi  \}  \big\}dx\\
&& +\int_{\partial M}\big\{\trv\{
12f_{e;m}\chi\psi\}+f_e\trv\{6\chi
\psi  _{ ;m}+6\psi_{;m}+6\chi  \gamma  _{ m} \gamma  _{
a} \psi  _{ ;a}\\
&&\qquad-12\chi  \psi L _{aa}-2\psi L _{ aa}-6\chi  \gamma  _{
m} \psi  \psi+3 \gamma  _{ m} \psi  \Psi  _{ T}\\
&&\qquad-3\chi  \gamma  _{
m} \psi  \chi  \psi +6 \chi  \gamma  _{
a}W _{ am} \} \big\}dy\bigg\}\,.
\end{eqnarray*}

Let $d:=\dim(V)$. We compute:\def\bork{\smallbreak\qquad}
\medbreak\bork$\trv\{6\psi_{;ii}-\tau\psi-6\gamma_i\gamma_jW_{ij}\psi-3\psi\psi\Psi\}=0$,
\bork$\trv(3\Psi_e\gamma_i\psi_o)_{;i}=12\trv(\tilde\gamma\gamma_m\tilde\gamma\gamma_m)(f_ef_o)_{;m}
    =12d(f_{e;m}f_o+f_{o;m}f_e)$,
\bork$\trv(3\Psi_o\gamma_i\psi_e)_{;i}=-6\trv(\tilde\gamma\gamma_m\gamma_m\tilde\gamma)(f_ef_o)_{;m}
     =6d(f_{e;m}f_o+f_{o;m}f_e)$,
\bork$\trv(6\psi_o\psi_{e;m}\gamma_m)=6\trv(\tilde\gamma\gamma_m\tilde\gamma\gamma_m)f_of_{e;m}
    =6df_{e;m}f_o$,
\bork$\trv(6\psi_e\psi_{o;m}\gamma_m)=6\trv(\tilde\gamma\tilde\gamma\gamma_m\gamma_m)f_ef_{o;m}
    =-6df_ef_{o;m}$,
\bork$f_e\trv\{(3\Psi\gamma_i\psi)_{;i}+6\psi\psi_{;}\gamma_i\}=12d(f_e^2f_0)_{;m}$.
\medbreak\noindent Thus interior integral is in divergence form.
Next we study the boundary term:
\medbreak\bork$f_e\trv\{6\psi_{;m}+6\chi\gamma_{m}\gamma_{a}\psi_{;a}
   -2\psi L _{ aa}-6\chi \gamma_{m}\psi\psi+6\chi\gamma_{a}W_{am}\}=0$,
\bork$f_{e;m}\trv\{12\chi\psi\}=-12df_{e;m}f_o$,
\bork$f_e\trv\{6\chi\psi_{ ;m}\}=-6df_ef_{o;m}$,
\bork$f_e\trv\{-12\chi  \psi L _{aa}\}=12df_ef_oL_{aa}$,
\bork$f_e\trv\{3 \gamma_m\psi\Psi_{
T}\}=3f_e^2f_o\trv\{\gamma_m\tilde\gamma\gamma_m\gamma_a\tilde\gamma\gamma_a
    +\gamma_m\tilde\gamma\gamma_a\tilde\gamma\gamma_m\gamma_a\}=18d f_e^2 f_o$,
\bork$f_e\trv\{-3\chi\gamma_{m}\psi\chi\psi\}
   =-3f_e^2f_o\trv\{\tilde\gamma\gamma_m\gamma_m\tilde\gamma\tilde\gamma\gamma_m\tilde\gamma\gamma_m
   +\tilde\gamma\gamma_m\gamma_m\tilde\gamma\gamma_m\tilde\gamma\gamma_m\tilde\gamma\}$
\newline\phantom{a}\qquad\qquad\qquad\qquad\qquad\qquad$
   =6df_e^2f_o$.
\medbreak\noindent Assertion (2) now follows after combining terms.
\end{proof}

\section{Orientability}\label{sect-6} We have assumed in previous
sections that $M$ is orientable. In fact, it is possible to relax
this restriction. We decompose $\partial M=N_1\cup...\cup N_\ell$
into connected components. Assume that each component $N_i$ is
orientable but that $M$ need not be orientable. On each component
$N_i$, choose an orientation to define $\tilde\gamma$ and thereby
define chiral bag boundary conditions. Let $\check M$ be the
orientable double cover; $\partial\check M=\check
N_1\cup...\cup\check N_{2\ell}$ where each $N_i$ lifts to two
distinct components in the double cover. Let $\check\gamma$ be
defined by an orientation of $\check M$.

Let $P$ be an operator of Dirac type on $M$. We can define $\tilde\gamma$
locally on $M$; we say that $P$ is {\it admissible} if
$\tilde\gamma P+P\tilde\gamma=0$ near any point $P$ of $M$; replacing
$\tilde\gamma$ by $-\tilde\gamma$ does not change this condition so
this is an invariantly defined notion. Equivalently, if $\check P$ is the
associated operator on $\check M$, this condition simply means
that $\check P\check\gamma+\check\gamma\check P=0$. Similarly, we say that
$F$ is ${\it admissible}$ if $F\tilde\gamma+\tilde\gamma F=0$
for any local orientation of $M$ or equivalently if
$\check F\check \gamma+\check\gamma\check F=0$. Since the heat trace invariants are
multiplicative under finite coverings, applying Theorem
\ref{thm-2.1} for $\{\check M,\check P,\check F\}$
(replacing if necessary $\chi_\theta$ by $-\chi_\theta$ on certain components if
necessary), then yields the following result for
$\{M,P,F\}$.
\begin{theorem}\label{thm-6.1}
Let $P$ be an operator of Dirac type on a compact non-orientable
smooth manifold of even dimension $m$ with orientable boundary. Let
$\BB_\theta$ be chiral bag boundary conditions defined
by orientations of the components of $\partial M$.  Let
$F\in C^\infty(\operatorname{End}(V))$. Assume that $P$ and $F$ are admissible. Then:
\begin{enumerate}
\item $\partial_\theta a_m(1,P,\BB_\theta)=0$.
\item $\partial_\theta a_m^\eta(1,P,\BB_\theta)=0$.
\item $\partial_\theta a_{m-1}^\eta(F,P,\BB_\theta)=0$.
\end{enumerate}
\end{theorem}

Similarly Theorem \ref{thm-3.1} extends immediately to yield:
\begin{theorem}\label{thm-6.2}
Let $P$ be an operator of Dirac type on a compact smooth
manifold of even dimension $m$ with orientable boundary $\partial M$. Let
$\BB_\theta$ be chiral bag boundary conditions defined by  orientations of the components of $\partial M$. Then
$a_m^\eta(1,P,\BB_\theta)=0$.
\end{theorem}

\section{Conclusions}\label{sect-7}
The crucial point of our arguments in Section \ref{sect-2} is that
we work with the associated first order operator. A chiral gauge
transformation is used to relate the original problem with
$\theta$-dependent boundary condition to a problem where the
operator is transformed to a $\theta$-dependent one and the
boundary conditions to a $\theta$-independent one. In a formula,
the chiral gauge transformation amounts to expressing
$P_{\theta,\BB_0}=e^{-\theta\tilde\gamma}P_{\BB_0}$, the domains
of the appropriate operators being crucial to understanding the
arguments involved. To emphasize: $\tilde\gamma
P_{\BB_0}\ne-P_{\BB_0}\tilde\gamma$. From the above stated
relation the variation of invariants with $\theta$ can be
evaluated and is given in Theorem \ref{thm-2.1}. As an application
of this theorem, we have given a very elegant proof of the
regularity of the eta invariant at $s=0$ in $m=4$ dimensions.
Furthermore, the coefficient $a_4 (1,\gamma_j \nabla_{e_j}
,\BB_\theta )$ in $m=4$ dimensions could be found using results
from standard mixed boundary conditions, as it does not depend on
$\theta$. Note, however, that in general $a_4 (1,P, \BB_\theta )$
shows non-trivial $\theta$-dependence, see Theorem \ref{thm-5.1}.
The eta invariant was computed \cite{G85,G87} for certain closed
even dimensional  pin${}^c$ manifolds; their $K$-theory groups and
some equivariant bordism groups were then determined. The results
of Section \ref{sect-6} show the eta invariant is well defined
with chiral bag boundary conditions even if $\partial
M\ne\emptyset$. However, the observation made in Remark
\ref{rmk-3.3} shows that the eta invariant is no longer a homotopy
invariant and thus a more careful study of the boundary
contribution needs to be made before the results cited above can
be extended to the case $\partial M\ne\emptyset$.

\section*{Acknowledgements} Research of
PG was partially supported by the MPI (Leipzig, Germany). KK
acknowledges support by the Baylor University Research Committee
and by the MPI (Leipzig, Germany).

\end{document}